\documentclass[a4paper,12pt]{article}
\usepackage[intlimits]{amsmath}
\usepackage[latin1]{inputenc}
\usepackage{amsfonts}
\usepackage{amssymb}
\usepackage{latexsym}
\usepackage{graphicx}
\usepackage{color}
\usepackage{graphics}
\usepackage{amsthm}
\title{Some remarks on depth of dead ends in groups}
\author{J\"org Lehnert}

\newtheorem{lemma}{Lemma}%[section]
\newtheorem{satz}[lemma]{Theorem}
\newtheorem{kor}[lemma]{Corollary}

\theoremstyle{definition}
\newtheorem{Def}[lemma]{Definition}

\newcommand{\ZZ}{\ensuremath{\mathbb{Z}}}
\newcommand{\NN}{\ensuremath{\mathbb{N}}}
\newcommand{\Hou}{\ensuremath{\mathrm{H}_2}}
\newcommand{\Houn}[1][n]{\ensuremath{\mathrm{H}_{#1}}}

\newcommand{\comment}[1]{}

\definecolor{dark}{gray}{.35}
\definecolor{light}{gray}{.7}
\begin{document}

\maketitle

\begin{abstract}
It is known, that the existence of dead ends (of arbitrary depth) in the Cayley graph of a group depends on the chosen set of generators. Nevertheless there exist many groups, which do not have dead ends of arbitrary depth with respect to any set of generators. Partial results in this direction were obtained by \v{S}uni\'c and by Warshall. We improve these results by showing that abelian groups only have finitely many dead ends and that groups with more than one end (in the sense of Hopf and Freudenthal) have only dead ends of bounded depth. Only few examples of groups with unbounded dead end depth are known. We show that the Houghton group \Hou\ with respect to a standard generating set is a further example. In addition we introduce a stronger notion of depth of a dead end, called strong depth. The Houghton group \Hou\ has unbounded strong depth with respect to the same standard generating set. 
\end{abstract}

\section{Introduction}
Let $G$ be a group and $X$ a finite set of generators. The (unoriented) Cayley graph $\Gamma=\Gamma(G,X)$ is the graph with vertex set $G$ whose edges are pairs $(g_1,g_2)\in G\times G$ with $g_1^{-1}g_2\in X^{\pm 1}$. Giving all edges the length $1$ we obtain a metric structure on $\Gamma$. We call this metric $d_X(\cdot,\cdot)$.

Many results on groups rely on the structure of geodesics in the Cayley graph. Because of the transitive action of G on the vertex set it suffices to consider geodesics from $1$ to each $g$. For some $g$ there might be no geodesic, that can be extended to a geodesic from $1$ to a $g'$ further away. Such elements $g$ are called \emph{dead ends} of $G$. More precisely: Let $n=d(1,g)$, $g$ is called a dead end of $G$ if the ball $B_g(1)$ of radius $1$ and center $g$ is contained in the ball $B_1(n)$ of radius $n$ with center $1$.

Let $k=\max\{l|B_g(l)\subseteq B_1(n)\}$, then $k$ is called the \emph{depth} of the dead end $g$. Dead ends and their depth were first considered by Bogopol'ski\v{i} in~\cite{bogo}, who proved that the depth of dead ends in a given hyperbolic group with a given set of generators is uniformly bounded. On the other hand~\cite{F} gives an example of a group with infinitely many dead ends, all of depth $2$, and \cite{lamp} gives an example of a group with arbitrary deep dead ends. 

It is easy to see that the property of having dead ends is not an invariant of a group. For example \ZZ\ generated by $\{2,3\}$ has the dead ends $1$ and $-1$. In fact \v{S}uni\'c~\cite{integers} proves that for each infinite group $G$ exists a generating set $X$, such that $G$ has dead ends with respect to $d_X$. Unfortunately even the property of having only dead ends of bounded depth is not a group invariant, as Riley and Warshall show in \cite{counter}\footnote{In this article Riley and Warshall give two examples of groups, with unbounded dead end depth with respect to one set of generators and bounded dead end depth with respect to some other set of generators. The proof of the finitely presented example is based on the uncorrected version of~\cite{finepres}. The authors have announced a corrected version.}.

But there are some results which do not depend on the set of generators. The result of Bogopol'ski\v{i} concerning hyperbolic groups was already mentioned above. In~\cite{integers} \v{S}uni\'c shows that \ZZ\ has only finitely many dead ends with respect to any generating set and Warshall~\cite{lattices} shows that for all weakly geodesically automatic groups, and hence all abelian groups, there exists a uniform bound on the depth of dead ends depending on the set of generators.

We will generalize the result of \v{S}uni\'c and a part of the result of Warshall by showing in Section~\ref{abelian}:
\begin{satz}\label{finite}
Let $G$ be an abelian group, generated by the finite set $X$. Then there exist only finitely many dead ends in $G$ with respect to $X$.
\end{satz}

The notion of ends of a graph goes back to Hopf~\cite{hopf} and Freudenthal~\cite{freuden}. Ends of a graph are equivalence classes of rays in the graph, where two rays are equivalent, if no finite set of vertices seperates them. In the case of Cayley graphs we consider w.l.o.g. only rays starting in $1$. The number of ends of a Cayley graph is known to be a quasi-isometry invariant. In particular we can speak of the number of ends of $G$. If $G$ is finite then $G$ has no ends. Hopf's result asserts that a finitely generated infinite group has one, two or infinitely many ends, and that two-ended groups are virtually $\ZZ$. Stallings celebrated structure Theorem describes the structure of groups with more than two ends in terms of amalgamated free products and HNN-extensions over finite amalgamated (associated) subgroups.

In Section~\ref{abelian} we will also show the following theorem.
\begin{satz}\label{ends}
Let $G$ be a finitely generated group with more than one end. Then there exists a uniform bound on the depth of dead ends depending only on the set of generators.
\end{satz}

There are not many examples of groups known, which have dead ends of arbitrary lage depth with respect to one set of generators. In Section~\ref{houghton} we prolong this list by showing that the Houghton group $\Hou$, which will be defined in Section~\ref{houghton} has arbitrary deep dead ends with respect to the standard generating set. We will also define the Houghton groups \Houn\ which we suspect to have arbitrary deep dead ends with respect to some sets of generators.

In the last Section we introduce the strong depth of a dead end. We show that \Hou\ has unbounded strong dead end depth.
 
\noindent{\it Acknowledgements:\hspace*{2mm}}
We would like to thank Robert Bieri for some usefull remarks.

\section{Dead ends in abelian groups and groups with more than one end}\label{abelian}
Let $G$ be a group and  $X$ a set of generators for $G$. Without loss of generality let $X$ be closed under inversion. Throughout this section we make frequent use of the canonical homomorphism $\pi$ which projects the free monoid $X^*$ over $X$ onto $G$ by sending each $x\in X$ to $x\in G$. We refer to the elements of the free monoid as words over the alphabet $X$. We call a word $w$, representing an element $g$ (i.e. $\pi(w)=g$), a geodesic word, if the length of $w$ (i.e. number of letters) equals $d_X(1,g)$. Then $w$ describes a geodesic path in $\Gamma$ from $1$ to $g$.

The main tool in the proof of Theorem~\ref{abelian} is the following observation: Assume that a word $w$ represents a dead end of a group. Then the word $w'\subset wX^*$ can only be geodesic if $w=w'$.

First of all we need the following observation concerning points in $\NN^n$.

\begin{Def}
Let $p=(p_1,p_2,\ldots,p_n),\ q=(q_1,q_2,\ldots,q_n)$ be points in $\NN^n$. We call $p$ and $q$ crossrelated ($p\lessgtr q$), if there exist $i,j$, such that $p_i<q_i$ and $p_j>q_j$. 
\end{Def}
\begin{lemma}\label{cross}
Any set $C$ of pairwise crossrelated points in $\NN^n$ has to be finite. 
\end{lemma}
\begin{proof}
We prove this by induction on $n$. Let $n=2$. Let $(x_1,y_1)\in C$. Then $C=C_1\bigsqcup C_2\bigsqcup\{(x_1,y_1\}$, where $C_1=\{(x_i,y_i)\in C|x_i<x_1,y_i>y_1\}$ and $C_2=\{(x_i,y_i)\in C|x_i>x_1,y_i<y_1\}$. $|C_1|<x_1$, because $(x_i,y_i),(x_j,y_j)\in C$ implies $x_i\neq x_j$, $|C_2|<y_1$, hence $|C|\leq x_1+y_1$.

Let $n=k+1$. Let $(d_1,d_2,\ldots d_n)\in C$. Define $C_j:=\{(c_1,c_2,\ldots,c_n)\in C|c_j<d_j\}$. Than $C=\bigcup_{j=1}^n C_j$. We only show that $C_1$ is a finite set, the rest follows analogously. For $0\leq l<c_1$ define \[D_l:=\{(c_2,c_3,\ldots,c_n)|(l,c_2,c_3,\ldots,c_n)\in C_1\}.\] $C$ is a set of pairwise crossrelated points, hence $C_1$ is a set of crossrelated points, hence all the $D_j$ are sets of pairwise crossrelated points in $\NN^k$. So $C_1$ is a finite disjoint union of finite sets and therefore finite.
\end{proof}

\begin{proof}[Proof of Theorem~\ref{finite}]
Let $G$ be an infinite abelian group generated by the set $X=\{x_1,x_2,\ldots,x_{2n}\}$, where $x_i=x_{n+i}^{-1}$. Let $w$ be a (geodesic) word representing an element $g\in G$. Then any permutation of the word $w$ is again a (geodesic) word representing the same element $g$. Hence it is sufficient to count the number of occurrences of the single letters and we can forget about the ordering of the letters. We obtain a surjective monoid homomorphism $\phi=\pi\circ\psi$ from $\NN^{2n}$ onto $G$ by defining the map (which is not a homomorphism) $\psi:\NN^{2n}\rightarrow X^*$,
\[\psi\left((i_i,i_2,\ldots,i_{2n})\right):=x_1^{i_1}x_2^{i_2}\ldots x_{2n}^{i_{2n}}.\]

Let $D=\{g_1,g_2,g_3, \ldots\}$ be the set of dead ends in $G$ and $w_j\in\phi^{-1}(g_j)$ such that $\psi(w_j)$ is a geodesic word. According to the observation at te beginning of this section, no $\pi$-preimage of a dead end can be a subword of the $\pi$-preimage ao another dead end. Therefore the $w_j$ have to be pairwise crossrelated and Lemma~\ref{cross} implies, that $D$ is a finite set.
\end{proof}

\begin{proof}[Proof of Theorem~\ref{ends}]
Let $X$ be set of generators for the group $G$, and $\Gamma=\Gamma(G,X)$ be a Cayley graph. Let $g\in G$ be a dead end. Because $G$ has more than one end, there exists a $k\in\NN$ such that $\Gamma\setminus B_1(k)$ has more than one component. By regularity of $\Gamma$, $\Gamma\setminus B_g(k)$ has more than one component, too. If $d(1,g)<k$ the depth of the dead end is $\leq 2k$. Otherwise $1$ lies in one component of $\Gamma\setminus B_g(k)$.
We choose a ray starting in $1$ whose end is in a different component. This ray hits the ball $B_g(k)$ at least once. Let $v_1$ denote the first intersection vertex of the ray and the ball and $v_2$ the last intersection vertex. Let $h$ be the vertex on the chosen ray and in the component of the end which has distance $k+1$ to $v_2$. Then $d(g,h)\leq 2k+1$ and 
\[d(1,h)\geq \stackrel{\leq d(1,v_1)}{\overbrace{(d(1,g)-k)}}+\stackrel{=d(v_2,h)}{\overbrace{k+1}}.\] 

Hence the depth of dead ends is bounded above by $2k$.
\end{proof}
\section{Dead ends in Houghton groups}\label{houghton}
The Houghton groups $\Houn$ were introduced by C.H.~Houghton in the late 70s (see~\cite{brown} for a more detailed discussion of Houghton groups). Although we only discuss $\Hou$ we define the whole family of groups.
\begin{Def}
Let $*^n=\{(k,l)|k\in\NN,l\in\{1,\ldots,n\}\}$ be the disjoint union of $n$ copies of \NN. The Houghton group \Houn\ is the group of all permutations $\phi$ of $*^n$ for which there exist constants $s(1),\ldots,s(n)\in\ZZ$ such that
\begin{eqnarray*}
\phi(k,l)&=&(k+s(l),l)
\end{eqnarray*}
for all but finitely many $(k,l)$.
\end{Def}
Then $\sum_{i=0}^n s(i)=0$ and \Houn[1] is just the symmetric group $S_\infty$, which is not finitely generated.
For $n>1$ all $H_n$ are finitely generated (in fact Brown showed in~\cite{brown} that $H_n$ is of type FP$_{n-1}$, but not of type FP$_n$). The shift $s_i\in H_n$ acts on $*^n$ as follows (rays modulo $n$):
\[
s_i(k,l)=\begin{cases}
         (k,l)    & i\neq l\neq i+1\\
         (k-1,l)  & l=i,k\geq 1\\
         (0,l+1)  & l=i,k=0\\ 
         (k+1,l)  & l=i+1
         \end{cases}
\]
For $n\geq 3$ the commutator $[s_i,s_{i+1}]$ acts as the transposition of $(0,i)$ and $(0,i+1)$. So the group generated by the shifts contains all finite permutations and is thereby equal to $\Houn$.

A special case is $\Hou$ which needs an additional generator, e.g. the transposition $\sigma=(0,1)\leftrightarrow(0,2)$. In this case $s_1=s_2^{-1}$ and we have in fact only one shift which we call $s$.
The group \Hou\ is isomorpic to $S_\infty\rtimes\ZZ$ and we can think of an element of $\Hou$ as a permutation with finite support followed by a shift.

According to this set of generators \Hou\ has the presentation
\begin{eqnarray}
\Hou &=&\left\langle s,\sigma\left|\sigma^2=(\sigma\sigma^s)^3=[\sigma,\sigma^{s^t}]=1,t\geq 2\right\rangle.\right.
\end{eqnarray}
Here $\sigma^s=s\sigma s^{-1}$ denotes the conjugation of $\sigma$ by $s$. The element $\sigma^{s^t}$ is the transposition of $(t,1)$ and $(t-1,1)$. 

In order to compute the length of a given group element with respect to this generating set we will use the following geometric interpretation of \Hou. Consider the biinfinite string of pearls numbered with nonzero integers, and a cursor which indicates the current position between two pearls. A word in \Hou\ is a sequence of movements of the cursor and commands to permute a finite number of pearls. The identity word is represented by all pearls in increasing order and the cursor between $-1$ and $1$. This cursor position is called origin. The element $s=s_1$ mentioned above moves the cursor one position to the right and the element $\sigma$ interchanges the pearls left and right to the cursor. One can easily check that the group elements of \Hou\ and these configurations are in 1-1 correspondence.

Let $X=\{s,\sigma\}$. As usual the distance $d_X(1,g)$ of a group element to the identity is defined as the minimum of the length of the words representing $g$ and a word is a sequence of commands to move the cursor one step or to interchange the pearls next to the cursor. The length of such a word equals the number of commands.

We introduce the following notation for elements of \Hou. As mentioned above each element is a permutation followed by shift. We write the permutation as a product of cycles and write the shift as index. For example $s=(1)_1$, $\sigma=(-1,1)_0$ or $g=(-1,2,4,-3,2)_{-3}$ (which means $g(-1)=s^{-3}(2)=-2$, $g(2)=1,\ g(4)=-6$ and so on).

The elements with trivial shift form a subgroup isomorphic to $S_\infty$ which is generated by $Y=\{\sigma_t=\sigma^{s^t}|t\in\ZZ\}$. Let $g_k=((-k,k)(-(k-1),k-1)\ldots(-3,3)(-2,2)(-1,1))_0$.
\begin{satz}\label{haupterg}
The element $g_k$ is a dead end of depth at least $k$ in \Hou\ with respect to the generating set $X$.
\end{satz}

We now have to calculate $d_X(1,g_k)$ and to find geodesic words representing $g_k$. We first calculate the distance $d_{Y_k}(1,g_k)$ in the subgroup $S_{2k}$ generated by $Y_k=\{\sigma_t|-k<t<k\}$. Let $M$ be the support of $g_k$, $M=\{-k\ldots -1,1,\ldots k\}$.
For each $g\in S_{2k}$ define the effect $e(g):=\sum_{i<j\in M}(\#\{(i,j)|g(i)>g(j)\})$. Than $e(g\cdot\sigma_t)=e(g)\pm 1$.
Hence $d_{Y_k}(1,g)\geq e(g)$. The effect of $g_k$ equals $e(g_k)=\sum_{i=1}^{2k}i-1=k(2k-1)$. If we find a word $w_k$ of length $k(2k-1)$ representing $g_k$ than $w_k$ is geodesic. 

The word $v_{l}:=\sigma_{-(l-1)}\sigma_{-(l-2)}\ldots\sigma_{l-2}$ represents the group element $(l-1,l-2,l-3\ldots,-l)_0$ and hence $u_l:=\sigma_{l-1}^{v_{l}}$ represents $(-l,l)_0$. The length of $u_l=1+4(l-1)$. $\pi(u_l)$ and $\pi(u_{l'})$ have disjoint support for $l\neq l'$ and therefore the concatenation of $u_l$ for all $0<l<k$ in any order represents $g_k$. Let $w_k=u_{k}u_{k-1}\ldots u_{1}$. Than the length of $w_k=\sum_{l=1}^{k}1+4(l-1)=k(2k-1)$. This word $w_k$ has a property which will be important later on. For each group element $g\in S_{2k}$ one can obtain a word representing $g$ by deleting some letters of $w_k$. From a Coxeter point of view this might not be surprising ($S_{2k}$ is a coxeter group) but can be easily seen without any of such arguments.
\begin{lemma}\label{loesch}
For any element $g\in S_{2k}$ one can obtain a word representing $g$ be deleting some letters in $w_k$.
\end{lemma}
\begin{proof}
We do this by induction on $k$. For $k=1$ there is nothing to show. Let $k=n+1$. We only need to show: It is possible to delete some of the letters of $u_1$ to obtain a word $x$ such that $\pi(x)(g^{-1}(k+1))=k+1$ and $\pi(x)(g^{-1}(-(k+1)))=-(k+1)$. This can be done by deleting some of the first letters each of $v_l$ and $v_l^{-1}$.
\end{proof}
We are now going to compute $d_X(1,g_k)$.  For each word $w$ in $Y$ (or in $Y_k$) let $w'$ be the word in $X$, which one obtains by replacing all $\sigma_t$ by $s^t\sigma s^{-t}$ and let $\tilde{w}$ be the reduced form of $w'$. Than $\tilde{u}_l=s^{-(l-1)}(\sigma s)^{2(l-1)}(\sigma s^{-1})^{2(l-1)}\sigma s^{l-1}$. 

How many commands are needed to interchange the pearl $-l$ and $l$ and bring the cursor back to the origin? W.l.o.g. let pearl $-l$ be moved before pearl $l$. The cursor has to come to the position $-(l-1)$ steps left of the origin for a first time which needs at least $l-1$ occurrences of $s^{-1}$. Afterwards at some stage the cursor has to come to position $l-1$ steps right to the origin which needs at least $2(l-1)$ occurrences of $s$. (Otherwise one of the elements would be a fixed point.) Finally the pearl $l$ has to be moved to its destination which needs again at least $2(l-1)$ occurrences of $s^{-1}$ and the cursor has to go back to the origin ($l-1$ times $s$). So we need at least $6(l-1)$ $s$ or $s^{-1}$. In addition we need at least (because of the effect) $4l-3$ occurrences of $\sigma$. Hence $\tilde{u}_l$ is a geodesic word.

The same argument shows that the only other geodesic word representing $(l,-l)_0$ can be obtained by interchanging $s$ and $s^{-1}$ in $\tilde{u}_l$.

\begin{lemma}\label{geod}
The word $\tilde{w}_k$ is a geodesic word representing $g_k$ and $d_X(1,g_k)=1+\sum_{l=2}^k (8l-5)$.
\end{lemma}
\begin{proof}
The length of the word $\tilde{w}=1+\sum_{l=2}^k (8l-5)$ ($\#\sigma's=\sum_{l=1}^k (4l-3)$, $\#s's=\sum_{l=2}^k (4l-2)$). We prove the statement by induction. For $k=1$ there is nothing to show. We only need to show the following claim: For any geodesic word $v_k$ representing $g_k$ there exist a word $v'_{k-1}$ representing $g_{k-1}$ and the $length(v'_{k-1})\leq length(v_k)-(8l-5)$.

As mentioned above we can think of $v_k$ as a sequence of commands to move the cursor or to swap the pearls next to the cursor. A geodesic word will never move the cursor outside of the region between pearl $-k$ and $k$. We enumerate the cursor positions inside by $c_{-(k-1)},\ldots,c_{(k-1)}$. So the origin is now called $c_0$. 

Let $v'_{k-1}$ be the word, one obtains by deleting the following letters of $v_{k}$: In the first step, we delete all letters $s$, which move the cursor away from or onto the origin in $v_{k}$. Afterwards we delete all letters $\sigma$ which interchange two pearls that were already interchanged before. 

How many letters $s$ are deleted in this procedure? For each pearl $l>1$ is moved through the origin which implies that the following situations have to occur: 1. The cursor is at position $c_1$ and the pearl $l$ is left next to the cursor and the cursor is moved to the origin. 2. The cursor is at the origin, the pearl $l$ is left next to the cursor and the cursor is moved to position $c_{-1}$. Hence for all pearls $l<1$ there are $2$ letters $s^{-1}$ deleted. Analogously for each pearl $l<-1$ the letter $s$ is deleted twice. In addition the cursor has to leave the origin a first time and come back to it a last time. All in all we have to delete at least $2(k-1)+2(k-1)+2=4-2$ letters $s^{\pm 1}$.

The word $v'_{k-1}$ represents $g_{k-1}$, because for $l>0$ one can check that $\pi(v'_{k-1})(l)=\pi(v)(l+1)+1$ and for $l<0$ that $\pi(v'_{k-1})(l)=\pi(v)(l-1)-1$. The number of $\sigma$s in $v'_{k-1}$ is (by definition) given by the effect of $\pi(v'_{k-1})$. Hence $length(v'_{k-1})\leq length(v_k)-(8l-5)$.
\end{proof}
\begin{proof}[Proof of Theorem~\ref{haupterg}]
We have to show that for all elements $g\in B_{g_k}(k)$ the distance $d(1,g)\leq d(1_{g_k})$. Let $\omega\in S_{\infty}$ and $t\in\ZZ$ such that $g=(\omega)_t$. Now $g\in B_{g_k}(k)$, hence $\omega\in S_{2k}$ and $|t|\leq k$. Lemma~\ref{loesch} implies that $d(1,g)\leq d(1,((k,-k)(k-1,-(k-1)\ldots(1,-1))_t))$. 
$g_k$ is an element of order $2$, so $g_k^{-1}=g_k$ and hence $\tilde{w}_k^{-1}$ ia another word representing $g_k$. This word ends with $k$ occurencies of $s$ and, as a consequence of symmetry there exist a geodesic word which represents $g_k$ and ends witk $k$ occurencies of $s^{-1}$. So in fact $d(1,((k,-k)(k-1,-(k-1)\ldots(1,-1))_t))\leq d(1,g_k)-|t|$.
\end{proof}
The situation in \Houn{} for $n>2$ is more complicated, but one still has a good combinatorial description of elements and generators and an argument similiar to the above argument for \Hou\ might show that \Houn{} (or at least \Houn[3]) has arbitrary deep dead ends.

\section{Strong depth of a dead end}
As mentioned above in \cite{counter} Riley and Warshall have shown that having dead ends of arbitrary depth is not a group invariant. One of their examples is the group with presentation
\[G=\langle a,t,u|a^2,[t,u],a^{-u}a^t;\forall i\in\ZZ, [a,a^{t^i}]\]
which has unbounded depth with respect to 
\[X=\{a,t,u,at,ta,ata,au,ua,aua\}\]
and depth bounded above by $2$ with respect to some different set of generators. Their examples of dead ends of depth $k$ have the following interesting property: Let $g$ be one of these dead ends and $n=d(1,g)$. There exists a geodesic from $g$ to an element $g'$ with distance $d(1,g')=n+1$ which never gets closer than $n-1$ to the identity. In other words: $g$ and $g'$ can be connected in $\Gamma(G,X)\setminus B_1(n-2)$. This yields to a new definition of depth of a dead end. To distinguish between them we will call the new one strong depth.

\begin{Def}
Let $\Gamma$ be the Cayley graph of a group $G$ with respect to a generating set $X$ and $g\in G$ a dead end with $d(1,g)=n$. The strong depth of $g$ is defined as the minimal number $k$ such that $g$ can be connected to a point of $\Gamma\setminus B_1(n)$ inside $\Gamma\setminus B_1(n-k)$. In other words: The strong depth of $g$ measures how many steps back to the identity a geodesic starting in $g$ has to take in order to leave the ball of radius $n$.
\end{Def}
The strong depth of a dead end is obviously less or equal the depth of it. The dead ends in \cite{counter} mentioned above are all of strong depth $2$. This is different from the situation in \Hou.
\begin{kor}
The element $g_k\in \Hou$ defined as in Theorem~\ref{haupterg} is a dead end of strong depth at least $k$. 
\end{kor}
\begin{proof}
In the proof of Theorem~\ref{haupterg} we have shown, that the distance of the form $(\omega)_t$, $\omega\in S_{2k}, |t|\leq k$ the distance $d(1,(\omega)_t)\leq d(1,g_k)$. Hence a geodesic from $g_k$ to a point outside the ball of radius $d(1,g_k)$ has to contain an element $h=(\omega')_{\pm k}$, $\omega\in S_{2k}$ and we have seen in proof of Theorem~\ref{haupterg} that $d(1,h)\leq d(1,g_k)-k$.
\end{proof}
To prove that the lamplighter group has unbounded strong depth only a few changes in the proof of~\cite{lamp} are needed. We do not know about any dead ends of large strong depth in the group of Riley and Warshall and therefore conclude by posing the following question: Does a group exist, which has unbounded strong depth of dead ends with respect to one set of generators but strong depth bounded above with respect to a different set of generators?


\begin{thebibliography}{HRRT}
\bibitem{bogo}{\it O.V. Bogopol'ski\v{\i}}, 'Infinite commensurable hyperbolic groups are bi-Lipschitz equivalent.', {\em Algebra and Logic} 36(3) (1997), 155-163.
\bibitem{brown}{\it K. Brown}, 'Finiteness properties of groups', {\em J. Pure and Applied Algebra} 44 (1987), 45-75.
\bibitem{F}{\it S. Cleary And J. Taback} 'Combinatorial properties of Thompson's group F.', {\em Transactions of the AMS} 356(7) (2004), 2825-2849.
\bibitem{lamp}{\it S. Cleary And J. Taback} 'Dead end words in lamplighter goups and other wreath products.', {\em The Quarterly Journal of Mathematics} 56(2) (2005), 165-178.
\bibitem{finepres} {\it S. Cleary And T. Riley}, 'A finitely presented group with unbounded dead end depth.', {\em Proceedings of the AMS}, 134(2) (2006), 343-349. Corrected version will be available on the Arxiv: {\it http://www.arxiv.org/abs/math.GR/0406443}.
\bibitem{freuden} {\it H. Freudenthal}, 'Neuaufbau der Endentheorie.', {\em Anals of Mathematics} 43 (1942), 261-279. 
\bibitem{hopf} {\it H. Hopf} 'Enden offener Räume und unedliche diskontinuierliche Gruppen.', {\em Comm. Math. Helv.} 15 (1943), 27-32.
\bibitem{counter}{\it T.Riley And A.D. Warshall}, 'The unbounded dead-end depth property is not a group invariant.', {\em International Journal of Algebra and Computation}, 16(5) (2006), 969-984. Corrected version in preperation.
\bibitem{integers} {\it Z. \v{S}uni\'{c}}, 'Frobenius problem and dead ends in integers', preprint.
\bibitem{lattices} {\it A.D. Warshall}, 'Deep pockets in Lattices and other groups', preprint.
\end{thebibliography}
\end{document}